\documentclass[11pt]{article}
\usepackage{amssymb}
\usepackage{amsbsy}
\usepackage{epsfig}
\usepackage{tikz}

\begin{document}
\pagestyle{myheadings}
\title{On the eccentricity energy of complete multipartite graph}
\author{
   {\em  Fernando Tura}\thanks{
        Partially supported by FAPERGS (Grant 17/2551-0000813-8).
        E-mail address:  \tt ftura@smail.ufsm.br} \\
    Universidade Federal de Santa Maria \\
    97105--900 Santa Maria, RS, Brazil \\
}


\def\floor#1{\left\lfloor{#1}\right\rfloor}

\newcommand{\diagonal}[8]{
\begin{array}{| c | r |}
b_i & d_i \\
\hline
 #1 & #5 \\
        & \\
 #2 & #6 \\
        & \\
 #3  & #7 \\
         &  \\
 #4  & #8 \\
 \hline
 \end{array}
}

\newcommand{\thresholdmatrix}[8]{
\left [
     \begin{array}{ccccccc}
            &   &   &   &  #1   &   & #2 \\
            &   &   &   &  \vdots   &   & \vdots \\
            &   &   &   &  #1  &   & #2 \\
     \\
         #3  & \ldots   & #3 &   &  #4    &   & #5 \\
     \\
         #6  & \ldots   & #6  &   &  #7   &   &  #8 \\
    \end{array}
\right ]
}

\newcommand{\lambdamin}{\lambda_{\min, n}}
\newcommand{\formulamin}{
 \frac{1}{2} (  (
   \lfloor \frac{n}{3} \rfloor \! - \! 1 ) \! - \! \sqrt{ ( \lfloor \frac{n}{3} \rfloor \! - \! 1)^2 \! +\!  4
   (n \! - \! \lfloor \frac{n}{3} \rfloor )
   \lfloor \frac{n}{3} \rfloor
  } )
}
\newcommand{\casei}{{\bf Case~1}}
\newcommand{\caseii}{{\bf Case~2}}
\newcommand{\caseiii}{{\bf Case~3}}
\newcommand{\caseiv}{{\bf Case~4}}
\newcommand{\subia}{{\bf subcase~1a}}
\newcommand{\subib}{{\bf subcase~1b}}
\newcommand{\subic}{{\bf subcase~1c}}
\newcommand{\subiia}{{\bf subcase~2a}}
\newcommand{\subiib}{{\bf subcase~2b}}
\newcommand{\subiic}{{\bf subcase~2c}}
\newcommand{\subiiia}{{\bf subcase~3a}}
\newcommand{\subiiib}{{\bf subcase~3b}}
\newcommand{\subiiic}{{\bf subcase~3c}}
\newcommand{\subiva}{{\bf subcase~4a}}
\newcommand{\subivb}{{\bf subcase~4b}}
\newcommand{\subivc}{{\bf subcase~4c}}
\newcommand{\myvar}{x}
\newcommand{\exvar}{\frac{\sqrt{3} + 1}{2}}
\newcommand{\Prf}{{\bf Proof: }}
\newcommand{\PrfSketch}{{\bf Proof (Sketch): }}
\newcommand{\boldQ}{\mbox{\bf Q}}
\newcommand{\boldR}{\mbox{\bf R}}
\newcommand{\boldZ}{\mbox{\bf Z}}
\newcommand{\boldc}{\mbox{\bf c}}
\newcommand{\sign}{\mbox{sign}}
\newcommand{\alphaseq}{{\pmb \alpha}_{G,\myvar}}
\newcommand{\alphaseqGprime}{{\pmb \alpha}_{G^\prime,\myvar}}
\newcommand{\alphaseqlam}{{\pmb \alpha}_{G,-\lambdamin}}
\newtheorem{Thr}{Theorem}
\newtheorem{Pro}{Proposition}
\newtheorem{Que}{Question}
\newtheorem{Con}{Conjecture}
\newtheorem{Cor}{Corollary}
\newtheorem{Lem}{Lemma}
\newtheorem{Fac}{Fact}
\newtheorem{Ex}{Example}
\newtheorem{Def}{Definition}
\newtheorem{Prop}{Proposition}
\def\floor#1{\left\lfloor{#1}\right\rfloor}

\newenvironment{my_enumerate}{
\begin{enumerate}
  \setlength{\baselineskip}{14pt}
  \setlength{\parskip}{0pt}
  \setlength{\parsep}{0pt}}{\end{enumerate}
}
\newenvironment{my_description}{
\begin{description}
  \setlength{\baselineskip}{14pt}
  \setlength{\parskip}{0pt}
  \setlength{\parsep}{0pt}}{\end{description}
}

\maketitle

\begin{abstract}
The eccentricity (anti-adjacency) matrix $\varepsilon(G)$ of a graph $G$ is obtained from the
distance matrix by retaining the eccentricities in each row and each column. The $\varepsilon$-eigenvalues of a graph $G$ are
those of its eccentricity matrix $\varepsilon(G),$ and the eccentricity energy (or the $\varepsilon$-energy) of $G$ is
the sum of the absolute values of $\varepsilon$-eigenvalues. In this paper, we establish some bounds for
the $\varepsilon$-energy of the complete multipartite graph $K_{n_1, n_2, \ldots,n_p}$  of order $n= \sum_{i=1}^p n_i $ and characterize the extreme graphs. 
 This partially answers the problem
given in Wang {\em et al.} (2019). We finish the paper showing graphs that are not $\varepsilon$-cospectral  with the same $\varepsilon$-energy.
\end{abstract}
{\bf keywords:} eccentricity matrix, eccentricity energy, complete multipartite graphs. 

\setlength{\baselineskip}{24pt}

\section{Introduction}
\label{intro}

Throughout this paper, all graphs are assumed to be finite, undirected and without loops or multiple edges.
Let $G = (V (G);E(G))$
be a graph with vertex set  $V (G) = \{v_1, v_2, \ldots, v_n\}$ and edge set $E(G).$
The distance $d_{G}(u, v)$ between two vertices $u$ and $v$ is the minimum length of the paths connecting them. By
convention, $d_{G}(v, v) = 0.$ Let $D(G) = (d_{uv})$ be the distance matrix of $G,$ where $d_{uv} = d_{G}(u, v).$ The eccentricity $\varepsilon(u)$ of the
vertex $u \in V(G)$ is given by $\varepsilon(u) = max\{d(u, v): v \in V(G)\}.$ The eccentricity matrix $\varepsilon(G) = (\epsilon_{uv} )$ of a graph $G$ is
defined as follows:

\begin{eqnarray}
\label{m(-1)}
 \epsilon_{uv}=  \left\{\begin{array}{ccccc}
          d_{G}(u,v), & if&  d_{G}(u,v) = min \{ \varepsilon(u), \varepsilon(v) \}&  &  \\
                 &     &       &    \\
          0,   &    otherwise.    &  &  & \\
   \end{array}\right.
 \end{eqnarray}  
which gives an equivalent definition of the $D_{max}$-matrix, due to Randi\'c \cite{Randic}.

The $\varepsilon$-matrix has nice application on Chemical Graph Theory.
For more details on this new matrix, we refer the
readers to \cite{Randic,Wang, Wang3}.

 Since $\varepsilon(G)$ is a symmetric matrix, then
eigenvalues of $\varepsilon(G)$ are real. Let $\varepsilon_1, \varepsilon_2, \ldots, \varepsilon_k$  be the distinct $\varepsilon$-eigenvalues.
Then the  $\varepsilon$-spectrum can be given as
$$spec_{\varepsilon}(G) = \{ \varepsilon_1^{m_1}, \varepsilon_2^{m_2}, \ldots, \varepsilon_k^{m_k} \}$$
where $m_i$ denotes the multiplicity of $\varepsilon$-eigenvalue $(i=1,2,\ldots, n).$

It is well-known that the graph energy is one important chemical index in Chemical Graph Theory. It is Gutman \cite{Gutman2012}
who firstly introduced the ordinary energy (or $A$-energy) of a graph $G$ as follows
\begin{equation}
 E_A(G) = \sum_{i=1}^n |  \lambda_i |
 \end{equation} 
where $\lambda_i$
is an $A$-eigenvalue of $G.$ Analogous graph energies invariants have been defined from the spectra of other graph
matrices \cite{Li2, Li7,Gutman,tura,steva} including the following eccentricity energy $E_{\varepsilon}(G)$ (or $\varepsilon$-energy) \cite{Wang2}:

\begin{equation}
\label{energy}
 E_{\varepsilon}(G)= \sum_{i=1}^n |  \varepsilon_i |.
 \end{equation}

A classical problem involving the energy is to find extremal graphs involving this parameter, that is to find
the graph on $n$ vertices with the largest energy. Concerning the eccentricity energy this is essentially an unheard problem. Based on this, Wang {\em et al.} \cite{Wang2} proposed the following problem: Given a set $\mathcal{S}$  of graphs, find an upper and lower bounds for the $\varepsilon$-energy of graphs in $\mathcal{S}$ and characterize
the extremal graphs.

To contribute to this problem, our main goal in this paper is given some bounds for
the $\varepsilon$-energy of the complete multipartite graph $K_{n_1, n_2, \ldots, n_p}$ on $n= \sum_{i=1}^{p} n_i$ vertices
and characterize the extreme graphs.

Since that finding not cospectral graphs with the same energy is an interesting problem in spectral graph theory, in this paper we continue this investigation  showing graphs that are not $\varepsilon$-cospectral  with the same $\varepsilon$-energy.


The paper is organized as follows. In Section 2 we describe some  known results about the  $\varepsilon$-spectrum of graphs.
In Section 3 we present the main results of the paper.
We finalize this paper, in Section 4, showing graphs with the same $\varepsilon$-energy.

\section{Notations and preliminaries}
\label{Diag}
As usual, let $K_{n-1,1}$ and $K_n$ denote the star and complete graph on $n$ vertices, respectively.
A complete split graph $CS(p_1,p_2)$ is a graph on $n=p_1+p_2$ vertices of a clique on $p_2$ vertices and 
and independent set on the remaining vertices in which each vertex of the clique is adjacent to each vertex of the independent set.

Let $M$ be a partitioned matrix

\begin{center}
$ M= \left[\begin{array}{cccc}
                M_{11}   &   M_{12}                       & \ldots    &  M_{1n}  \\
               M_{21}    &   M_{22}                           & \ldots           & M_{2n}\\
            \ldots         &         \ldots            &  \ldots       & \ldots     \\
        M_{n1}       &    M_{n2}   &\ldots         &     M_{nn} \\
\end{array}\right]$
\end{center}
where each $M_{ij}$ is a submatrix (block)  of $M.$ If $q_{ij}$  denotes  the average  row sum  of $M_{ij}$  then the matrix  $\mathcal{Q}= (q_{ij})$ is called a {\em quotient matrix} of $M.$  If the row sum of each  block $M_{ij}$ is a constant then the partition is called  {\em equitable}.  The following is a well known result on equitable partition of matrices (see, for example, Lemma 2.3.1, \cite{Bro}).

\begin{Lem}
\label{lem1}  Let $\mathcal{Q}$  be a quotient matrix of a square matrix $M$ corresponding to an equitable partition. Then the spectrum of $M$ contains  the spectrum of $\mathcal{Q}.$ 
\end{Lem}

The next result  shows a relation between the  $A$-eigenvalues and $\varepsilon$-eigenvalues, according \cite{Wang}.

\begin{Lem}
\label{lem2}
Let  $G$ be a connected graph
of order $n,$ diameter $d = 2$ and maximum degree  $\Delta < n -1.$ 
Then \begin{equation}
\label{eq5}
\varepsilon(G)= 2 A(\overline{G})
\end{equation}
where $\overline{G}$ is the complement of $G.$
\end{Lem}

We finalize this section presented a result that will be important in the sequence of the article due \cite{Lin}.

\begin{Lem}
\label{lem3}
The graph $K_{n-1,1}$ is the unique graph, which have maximum
distance spectral radius $\rho(D(K_{n-1,1}))$ among all graphs with diameter 2.
\end{Lem}

\section{Main Results}
\label{Diag}

Let $n= \sum_{i=1}^n n_i$ with $n_1 \geq n_2 \geq \ldots \geq n_p \geq 1$ and $K_{n_1, n_2, \ldots,n_p}$ denotes the complete multipartite graph.
In this section, for fixed $n,$ we obtain some bounds for
the $\varepsilon$-energy of the complete multipartite graph
and characterize the extreme graphs.

\begin{Thr}
\label{main1}
Let $G= K_{n_1, n_2, \ldots,n_p}$ be the multipartite graph of order $n\geq 4$ with $n_1\geq n_2 \geq \ldots \geq n_p \geq 1.$ 
Then the $\varepsilon$-eigenvalues of $K_{n_1, n_2, \ldots,n_p}$ are:
\begin{my_description}
\item[i.]  $n_i -1$ and $-2$
 with multiplicity $n_i -1,$ 
if $n_i \geq  2$ $(i=1,2,\ldots,p).$

\item[ii.] $n-1$ and $-1$ with multiplicity $n-1,$ if $n_i =1$ $(i=1,2,\ldots,p).$

\item[iii.] $-2$ with multiplicity $p_1 -1,$ $-1$ with multiplicity $p_2 -1$ and the roots $\varepsilon_1, \varepsilon_2$ of equation
\begin{equation}
\label{eq5}
\varepsilon^2 -(2p_1 + p_2 -3)\varepsilon +p_1p_2 -2p_1-2p_2+2
\end{equation}
where $p_1 = \sum_{i=1}^j n_i$ and $p_2= \sum_{i=j+1}^p n_i,$ if $n_1 \geq n_2 \geq \ldots \geq n_j \geq 2$ and $n_{j+1}= \ldots =n_p=1.$ 
\end{my_description}
\end{Thr}

\noindent{\bf Proof:}
$i.$ If $G$ has $n_i \geq 2$ $(i=1,2,\ldots,p),$ by Lemma \ref{lem2} the $\varepsilon$-eigenvalues are obtained from the eigenvalues of $\overline{G}.$
Since the  complementar is the disjoint union of some complete graphs $K_{n_i}$ $(i=1,2,\ldots,p)$ then follows the item $i.$

$ii.$ If $n_i =1$ $(i=1,2,\ldots,p)$ then the complete multipartite graph is isomorphic to the complete graph $K_n.$
In this case, the $\varepsilon$-spectrum coincides with the $A$-spectrum and the result follows.

$iii.$ Now we assume that the complete multipartite graph has vertices  $n_1 \geq n_2 \geq \ldots \geq n_j \geq 2$ and $n_{j+1}= \ldots =n_p=1.$ 
In this case, the graph is isomorphic to the complete split graph $CS(p_1, p_2)$  where  $p_1 = \sum_{i=1}^j n_i$ is the set of independent vertices
 and $p_2= \sum_{i=j+1}^p n_i$ is the set of clique vertices. We have that
the eccentricity matrix of $G$ is
\begin{center}
$ \varepsilon(G)= \left[\begin{array}{ccc}
                2 (J - I)_{p_1 \times p_1}    &                       & J_{p_1 \times p_2}    \\
               J_{p_2 \times p_1}    &                          & (J -I )   \\
\end{array}\right].$
\end{center}

Since the independent set of vertices contributes for the $\varepsilon$-eigenvalue $-2$ while the clique set of vertices contributes for the $\varepsilon$-eigenvalue $-1$
( see,  Theorem 2.34 of \cite{distance}) and
 the given partition of the matrix  is equitable, hence from the quotient matrix,
we get the characteristic polynomial  is the following
$$p(\varepsilon) =(\varepsilon+2)^{p_1} (\varepsilon+1)^{p_2} q(\varepsilon)$$
where $q(\varepsilon)$ is given by equation (\ref{eq5}), and the result follows. $\hspace{2,5cm} \square$



\begin{Thr}
\label{main2}
Let $G= K_{n_1, n_2, \ldots,n_p}$ be the multipartite graph of order $n\geq 4$ and $\varepsilon_1(G)$ denotes the $\varepsilon$-spectral radius.
Then
\begin{equation} 
\varepsilon_1(G) \leq  (n -2) +\sqrt{n^2-3n+3}
\end{equation}
 and equality holds if and only if $G\cong K_{n-1,1}.$
\end{Thr}
\noindent{\bf Proof:} Let $K_{n-1,1}$ be the star on $n$ vertices with $\varepsilon$-spectral radius equals $(n -2) +\sqrt{n^2-3n+3},$
obtained from equation (\ref{eq5}), and let $G$ be the complete multipartite graph with maximum $\varepsilon_1(G).$
Since in this case, the eccentricity matrix and the distance matrix are the same, by Lemma \ref{lem3}, $\varepsilon_1(G)= \rho(D(G)) \leq \rho(D(K_{n-1,1})) = (n -2) +\sqrt{n^2-3n+3}. \hspace{3,5cm} \square$

Before we provide an upper and lower bounds for the $\varepsilon$-energy of complete multipartite graph,
we need of the following technical result.

\begin{Lem}
\label{lem4}
Let $b$ and $c$ be positive real numbers, and let $\varepsilon^2 -b\varepsilon +c$ have real roots $\varepsilon_1,\varepsilon_2.$
Then $| \varepsilon_1 | + | \varepsilon_2| = b.$
\end{Lem}
\noindent{\bf Proof:}  Since $b > \sqrt{b^2 -4c},$ both $\varepsilon_1= \frac{b +\sqrt{b^2 -4c}}{2}$ and $\varepsilon_2= \frac{b -\sqrt{b^2 -4c}}{2}$ are positive.
Thus $| \varepsilon_1 | + | \varepsilon_2| =  \varepsilon_1  +  \varepsilon_2  =  b. \hspace{6,2cm} \square$

\begin{Thr}
\label{main3}
Let $G= K_{n_1, n_2, \ldots,n_p}$ be the multipartite graph of order $n\geq 4$ with $n_1\geq n_2 \geq \ldots \geq n_p \geq 1.$ 
Then the $\varepsilon$-energy of $G$ satisfies:
 \begin{equation}
 \label{eq6}
2n -2 \leq E_{\varepsilon} (G) \leq 2(n -2) +2\sqrt{n^2-3n+3}
\end{equation}

\end{Thr}

\noindent{\bf Proof:}  
If $G\cong K_n,$ by Theorem \ref{main1}, we have $E_{\varepsilon}(G)= 2(n-1),$
that is, the $\varepsilon$-energy of $G$ satisfies the inequality  (\ref{eq6}).

If $G= K_{n_1, n_2, \ldots,n_p}$ has $n_i\geq 2$ $(i=1,2,\ldots,p),$ by Theorem \ref{main1},
$\varepsilon$-eigenvalues of $G$  are  $2(n_1 -1), 2(n_2 -1), \ldots, 2(n_p -1), -2^n,$ then 
\begin{equation}
\label{eq} 
 E_{\varepsilon}(G) = 4(n_1 +n_2 +\ldots+n_{p} -p) = 4(n-p)
 \end{equation}
  Then the $\varepsilon$-energy decreases when $p$ increases.
Hence $K_{2,2,\ldots,2}$  has the minimum $\varepsilon$-energy which is equal to $2n,$ while that $K_{n_1, n_2}$ has the maximum $\varepsilon$-energy equals to $4(n-2).$
The $\varepsilon$-energy of both graphs satisfy the inequality on left in (\ref{eq6}).
For showing the inequality on the right of (\ref{eq6}), note 
$$E_{\varepsilon}(K_{n_1,n_2}) = 4(n-2) = 2( n-2 +n-2) \leq 2( n-2)  + 2\sqrt{n^2 -3n+3}.$$

Now, we assume that $G\cong CS(p_1, p_2)$ with $n=p_1 +p_2$  where $p_1 = \sum_{i=1}^j n_i$  and $p_2= \sum_{i=j+1}^p n_i.$ 
 By Theorem \ref{main1}, the $spec_{\varepsilon}(G) = \{ \varepsilon_1, \varepsilon_2, -2^{p_1-1}, -1^{p_2-1} \},$
where  $\varepsilon_1$ and $\varepsilon_2$ are the roots of equation (\ref{eq5}).
We assume that $G \neq K_{n-1,1},$ that is  $p_1, p_2 \geq 2.$  Then consider two cases:

{\em Case 1:}  $\varepsilon_1 \varepsilon_2 < 0.$  
Assuming that $\varepsilon_2 < 0 $ (or $-\varepsilon_2 >0$)  follows that $E_{\varepsilon}(G) = 2 \varepsilon_1.$
For the inequality on left in $(\ref{eq6}),$ since that $\varepsilon_1 + \varepsilon_2 = 2(p_1 -1) + p_2 -1= n-1 +p_1 -2$ then
$$E_{\varepsilon}(G) =  2(n-1 +p_1 -2 -\varepsilon_2)  > 2(n-1).$$
For the inequality on right in $(\ref{eq6}),$ by Theorem \ref{main2} follows 
$$E_{\varepsilon}(G) =  2 \varepsilon_1 < 2(n -2) +2\sqrt{n^2-3n+3}.$$

{\em Case 2:}  $\varepsilon_1 \varepsilon_2  \geq 0.$ By Lemma \ref{lem4}, we have $|\varepsilon_1| +| \varepsilon_2| = 2p_1 + p_2 -3.$
Then the $\varepsilon$-energy of $G$ satisfies $E_{\varepsilon}(G) = 2 ( 2p_1 +p_2 -3) = 2 (n-1 +p_1 -2) \geq 2(n-1),$ with equality if $p_1=2.$
For the inequality on right in $(\ref{eq6}),$ by Theorem \ref{main2} again follows 
$$E_{\varepsilon}(G) =   \varepsilon_1 + \varepsilon_2 < 2(n -2) +2\sqrt{n^2-3n+3}.$$

If $G\cong K_{n-1,1}$ since $\varepsilon$-eigenvalues of $G$  are $(n-2) \pm \sqrt{n^2 -3n+3}, -2^{n-2}$ then $G$ is the graph that reaches the upper bound in the
inequality (\ref{eq6}) and the result follows.
$\hspace{9,3cm} \square$

\begin{Cor}
Among all connected complete multipartite graphs $K_{n_1, n_2, \ldots,n_p}$ on $n$ vertices,
the complete graph $K_n$ and the complete split graph $CS(2, n-2)$ have the minimum $\varepsilon$-energy. Furthermore,
the star $K_{n-1,1}$ has the maximum $\varepsilon$-energy.
\end{Cor}

\section{ $\varepsilon$-equienergetic graphs}

Let $G$ be a connected graph with $u\in V(G).$ The set of vertices at distance $r$ from $u$ is denoted by $G_{r}(u)$ $(0\leq r\leq \varepsilon(u)).$
A graph $G$ of diameter $d$ is called {\em antipodal} if, for any given vertex $u\in V(G),$ the set $\{u \} \cup G_{d}(u)$ consists of vertices which are mutually at distance $d.$
In other words, there exists a partition of the vertex set into classes (called the fibres of $G$) with the property that two
distinct vertices are in the same class iff they are at distance $d.$ If all the fibres have the same cardinality, say $a,$ we say that $G$ is an $a$-antipodal graph \cite{Wang3}.
Obviously, the regular complete multipartite graph $K_{a,a,\ldots,a}$ is an $a$-antipodal graph.

Recall that the strong product $G\boxtimes H$ of two graphs  $G=(V,E)$ and $H=(W,F)$ is the graph with set vertex  $V\times W$ for which  $(v_1,w_1)$ and $(v_2,w_2)$ are adjacent if and only if $v_1=v_2$  and $w_1=w_2$ or $v_i w_i \in E(G_i)$ for $i=1,2,$ according \cite{Hammack}.

The following is a known result about  the $\varepsilon$-spectrum of the strong product (see, \cite{Wang2}).

\begin{Thr}
\label{main5}
 Let $G$ be an $a$-antipodal graph with order $m,$ diameter $d$ and $\varepsilon$-eigenvalues $d(a-1), -d$ with multiplicity  $ \frac{m}{a}, \frac{m(a-1)}{a},$ respectively, where $a|m.$
 Let $H$ be a connected graph with order $n$ and $diam(H) < d.$ Then the $\varepsilon$-eigenvalues of $G\boxtimes H$ are $dn(a-1), 0,-dn$ with multiplicity $\frac{m}{a}, m(n-1), \frac{m}{a}(a-1),$ respectively.
 \end{Thr}

\begin{Thr}
\label{main6}
Let $K_{n,n,n,n}$ and $K_{n,n}$ be two $n$-antipodal graphs of order $4n$ and $2n,$ respectively.
Then $K_{n,n} \boxtimes K_2$ and $K_{n,n,n,n}$ are not $\varepsilon$-cospectral and $\varepsilon$-equienergetic graphs.
\end{Thr}
\begin{Prf} 
Obviously, $K_{n,n} \boxtimes K_2$ and $K_{n,n,n,n}$ are not $\varepsilon$-cospectral, since that $0 \in Spec_{\varepsilon}(K_{n,n} \boxtimes K_2)$
while that  $0 \notin Spec_{\varepsilon}(K_{n,n,n,n}).$
By direct computation, we have that $E_{\varepsilon}(k_{n,n}\boxtimes K_2)= 2mnd(\frac{a-1}{a})= 2(2n)(4)\frac{n-1}{n}=16(n-1)=4(4n-4)=E_{\varepsilon}(K_{n,n,n,n}). \hspace{8cm}\square$
\end{Prf}

Note that from equation (\ref{eq}), it is easy to see that two non isomorphic complete multipartite graphs of order $n$
with same parity $p$ are not $\varepsilon$-cospectral and $\varepsilon$-equienergetic graphs. Using the Theorem \ref{main6} follows.

\begin{Cor}
The graphs $K_{n,n} \boxtimes K_2$ and $K_{n+i,n,n,n-i}$ are not $\varepsilon$-cospectral and $\varepsilon$-equienergetic graphs,
for $i=1,2, \ldots,n-2$ and $n\geq 4.$ 
\end{Cor}

\setlength{\baselineskip}{14pt}

\end{document}